\documentclass{amsart}

\usepackage{amssymb} 
\usepackage{epsf}

\newcommand{\bbN}{\mathbb{N}} \newcommand{\cW}{\mathcal{W}}
\newcommand{\cV}{\mathcal{V}} \newcommand{\cS}{\mathcal{S}}
 \newcommand{\cP}{\mathcal{P}}
\newcommand{\cM}{\mathcal{M}} 

\newcommand{\eps}{\varepsilon} \newtheorem{lemma}{Lemma}
\newtheorem{theorem}{Theorem}

\newcommand{\df}{\operatorname{def}}
\begin{document}

\title[Note on P-resolutions]{A note on P-resolutions of cyclic
  quotient singularities}

\author{Ludwig Balke}
\address{Math.~Institut University Bonn\\ Beringstr.~1\\53115 Bonn\\
  GERMANY}

\email{Balke@math.uni-bonn.de}


\begin{abstract}
  P-resolutions of a cyclic quotient singularity are known to be in
  one-to-one correspondence with the components of the base space of
  its semi-universal deformation \cite{Kollar-Shepherd-Barron}.
  Stevens \cite{Stevens} and Christophersen \cite{Christophersen} has
  shown that P-resolutions are parametrized by so-called chains
  representing zero or, equivalently, certain subdivisions of polygons.
  I give here a purely combinatorial proof of the correspondence
  between subdivisions of polygons and P-resolutions.
\end{abstract}

\maketitle

\section{Introduction}
The discovery of Koll\'ar and Shepherd--Barron in
\cite{Kollar-Shepherd-Barron} that the reduced components of the
versal base space of a quotient singularity are in one-to-one
correspondence with the P-resolutions of this singularity provided a
more conceptual understanding of the deformation theory of cyclic
quotient singularities. Christophersen \cite{Christophersen} and
Stevens \cite{Stevens} found a beautiful description of these
components in terms of so-called chains respresenting zero or
subdivisions of polygons. Nevertheless, their proofs were very
involved and used algebro-geometric methods. Altmann \cite{Altmann}
gave a description of P-resolutions and their correpondence to chains
representing zero in terms of toric varieties.

In this note, I will emphasize that P-resolutions on the one hand
(cf.~Lemma~\ref{lem_Presolution}) and subdivisions of polygons on the
other hand are purely combinatorial objects. Hence it is interesting
to establish the correspondence between both classes by purely
combinatorial methods. Therefore, I will refer to cyclic quotient
singularities and P-resolutions only in so far as it is necessary to
define the objects to be investigated here. The reader can find in
\cite{Stevens} more information on cyclic quotient singularities and
their deformation theory.

\section{Sequences and subdivisions}

Let $\cW$ denote the monoid of sequences $(a_1,\dots,a_k)$ with
$a_i\in \bbN\setminus\{0\}$ and $\cW_2$ the submonoid of such
sequences with $a_i\geq 2$. The empty sequence $()$ is denoted by
$\eps$.

Furthermore, let $M$ denote the free monoid with generators
$\alpha,\beta$ and $M_2$ the free monoid with generators $\alpha,
\gamma$. Both monoids act on $\cW$ respectively $\cW_2$ in the
following manner from the lefthand side or from the righthand side,
respectively :

\begin{alignat*}{2}
  \alpha \eps & = \eps &  \eps\alpha & = \eps\\
  \alpha(a_1,\dots,a_k)& =(a_1,\dots,a_k +1) &
  \qquad (a_1,\dots,a_k)\alpha& =(a_1+1,\dots,a_k) \\
  \beta a & = a(1) & a\beta & = (1)a\\
  \gamma a & = \alpha \beta a = a(2) & a\gamma & = a \beta\alpha =
  (2)a
\end{alignat*}

Here $a$ is an arbitrary element of $\cW$ respectively $\cW_2$ and the
composition in the respective monoids is denoted by juxtaposition. The
following lemma is easy to observe:

\begin{lemma}
  For each $a\in\cW$ there exists a unique element $\rho\in M$ such
  that $a=\rho \eps$ and $\rho=1$ or $\rho=\rho'\beta$ for a suitable
  $\rho'$. An analogous statement holds for $\cW_2$ and $M_2$.
  
  The inversion of a sequence in $\cW$ is defined by
  $\overleftarrow{(a_1,\dots,a_k)}:=(a_k,\dots,a_1)$. Analogous
  inversions are defined for the free monoids $M$ and $M_2$. Then we
  have for $a\in\cW, \rho \in M$ or $a\in\cW_2,\rho\in M_2$,
  respectively.
\[ \rho a = \overleftarrow{a}\overleftarrow{\rho}\]
\end{lemma}

We define a map $R:\cW\to\cW$ inductively as follows.  On $\cW_2$ we
have $R(\eps)=\eps$ and
\[R(\alpha a) =\gamma R(a)\quad\text{and}\quad
R(\gamma a) = \alpha R(a)\] for $a\in\cW_2, a\neq\eps$. For arbitrary
elements $a\in\cW$, the map $R$ is defined by the rule
\[ R(a(1)a') = R(a)(1)R(a').\]

A straightforward induction shows:

\begin{lemma}
  $R$ is an involution, i.e.~$R^2=\mbox{id}_{\cW}$, and
\[
R(\overleftarrow{a})=\overleftarrow{R(a)}.
\]\label{lem_involution}
\end{lemma}

\begin{lemma}
  For $a',a''\in\cW$ and $e,f\geq 2$:
\[ R(a' (e+f-1) a'') = R(a'(e))R ((f)a'').\]
\label{lem_split}
\end{lemma}

A nonempty sequence $(a_1,\dots,a_k)\in\cW_2$ defines a continued
fraction
\[ [a_1,\dots,a_k] =a_1 - \cfrac{1}{a_2 -
  \cfrac{1}{\dots - \cfrac{1}{a_k }}}
\]
It is an easy consequence from Lemma~\ref{lem_involution} that
$[R(a)]=\frac{n}{n-q}$ if $[a] =\frac{n}{q}$. Hence the operator $R$
generalizes Riemenschneider's point diagram rule to $\cW$
(cf.~\cite{Riemenschneider}, \cite{Stevens}, 1.2).

For later use, we need a relation between the positions in the
sequence $a$ and the positions in the sequence $R(a)$. This relation
will be inductively defined. Let $\ell$ be the length of $a$ and $r$
the length of $R(a)$. The $(\ell+1)$-th position of $\gamma a$ is
associated with the $r$-th position of $R(\gamma a)=\alpha R(a)$ and
the $\ell$-th postion of $\alpha a$ is associated with the $(r+1)$-th
position of $R(\alpha a) = \gamma R(a)$.  Furthermore, the $(\ell
+1)$-th position of $\beta a$ is associated with the $(r+1)$-th
position of $R(\beta a)=\beta R(a)$.

Obviously, the $i$-th position of $a$ is associated with the $j$-th
position of $R(a)$ if and only if the $j$-th position of $R(a)$ is
associated with the $i$-th position of $R(R(a))=a$.

We fix a countable set $\cV$ together with a linear ordering $\leq$ of
$V$.  Elements of $V$ will be called vertices and if $V$ is a finite
subset of $r$ vertices then $v_i$ with $1\leq i\leq r$ is defined by
the equation $V=\{v_0< \dots <v_r\}$.

A {\em subdivision of an polygon with distinguished vertex}, or
shortly subdivision, is a pair $(V,\Delta)$ such that $V$ is a
nonempty finite subset of $\cV$ and $\Delta$ is a set of subsets of
$V$ with the following property.  There exists a subdivision of a
plane polygon into triangles and a bijection between the vertices of
the polygon and the elements of $V$ such that $v_i$ and $v_{i+1}$
correspond to consecutive vertices and each set of vertices of a
triangle is identified with an element of $\Delta$ and vice versa. We
accept a digon as plane polygon, hence $(V,\emptyset)$ is a
subdivision if $\# V =2$. An {\em isomorphism} between subdivisions
$(V,\Delta)$ and $(V',\Delta')$ is a bijection from $V$ to $V'$ which
maps $\Delta$ bijectively onto $\Delta'$. Hence this bijection either
preserves or reverses the ordering $\leq$.  Obviously, the notion of
subdivision can be defined in a purely combinatorial way without
reference to plane polygons. But we do not bother the reader with an
axiomatic approach since this would not enhace the understanding of
this concept.

We define two laws of compositions on the set of subdivisions.  Let
$(V_1,\Delta_1)$ and $(V_2,\Delta_2)$ be subdivisions with
\begin{eqnarray*}
v^0&:=&\min V_1 = \min V_2\\
v^1&:=&\max V_1 < v^2:=\min V_2\setminus \{\min
V_2\}.
\end{eqnarray*}
Then we can define a new subdivision
$(V_1,\Delta_1)(V_2,\Delta_2):=(V,\Delta)$ with
\[V:=V_1\cup V_2\quad\text{and}\quad \Delta := \Delta_1\cup
\Delta_2\cup \{\{v^0,v^1,v^2\}\}\] A subdivision is called {\em
  irreducible} if it is not the product of two subdivisions and it is
called {\em primary} if there is atmost one triangle with $\min V$ as
vertex. Obviously, a primary subdivision is also a irreducible one,
but, in general, the converse is not true.  Furthermore, each
subdivision has a unique decomposition into irreducible factors.

Assume now $v^1=v^2$, then $(V_1,\Delta_1)*(V_2,\Delta_2):=(V,\Delta)$
with
\[V:= V_1\cup V_2\quad\text{and}\quad \Delta:=\Delta_1\cup\Delta_2\]
Each irreducible subdivision can by uniquley decomposed into primary
factors with respect to the law of composition $*$.

Let $(V,\Delta)$ be a subdivison of an $n$-gon and
$a=(a_1,\dots,a_n)\in\cW$. The sequence $a$ can be understood as a
labelling of the vertices $v$ of $V$ and we define the {\em degree}
$\deg v=\deg_a v$ of a vertex with respect to this labelling in such a
way that it emphasizes the special role of vertices with label $1$.
Consider triangles of $\Delta$ as equivalent if they have an edge in
common and one of the triangles has a vertex with label $1$. Then the
degree of a vertex $v$ is the number of equivalence classes of
triangles having $v$ as one of its vertices. In particular, a vertex
with label $1$ has always degree $1$.  We call $a$ {\em admissible}
for $(V,\Delta)$, if we always have
\[\deg_a v_i \leq a_i.\]
The difference $\df_a v_i := a_i - \deg_a v_i$ is called the {\em
  defect} of the vertex $v_i$.  The set of triples $(V,\Delta,a)$ with
$(V,\Delta)$ a subdivision and $a$ admissible for $(V,\Delta)$ is the
set of {\em labelled subdivisions} and is denoted by $\cS$.

A vertex $v$ of $(V,\Delta)$ is called a {\em splitting vertex} if it
is the maximal or the second vertex of an irreducible factor of
$(V,\Delta)$. A vertex $v$ is called {\em saturated} if its defect is
$0$, or if it is a splitting vertex and both neighbouring vertices
have a label different from $1$, or if this vertex together with its
both neighbours forms a triangle.  A triangle $t\in\Delta$ is called an
{\em interior triangle} if and only if $t=\{v_i,v_j,v_k\}$ with
$i+1<j<k-1$ and $0<i$.

A nonempty sequence $a\in\cW$ can be associated with a labelled graph
$\Gamma (a)$ as follows. If $a=(a_1,\dots,a_k)$ then $\Gamma (a)$ is a
chain with $k$ nodes which are labelled consecutively with the numbers
$-a_1,\dots,-a_k$. For $a\in\cW_2$ with $[a]=\frac{n}{q}$, the graph
$\Gamma (R(a))$ is the graph of the minimal resolution of the cyclic
quotient singularity $X(n,q)$ (cf. \cite{Riemenschneider}).

\begin{figure}

    \epsfysize3cm
    \hspace*{\fill}\epsffile{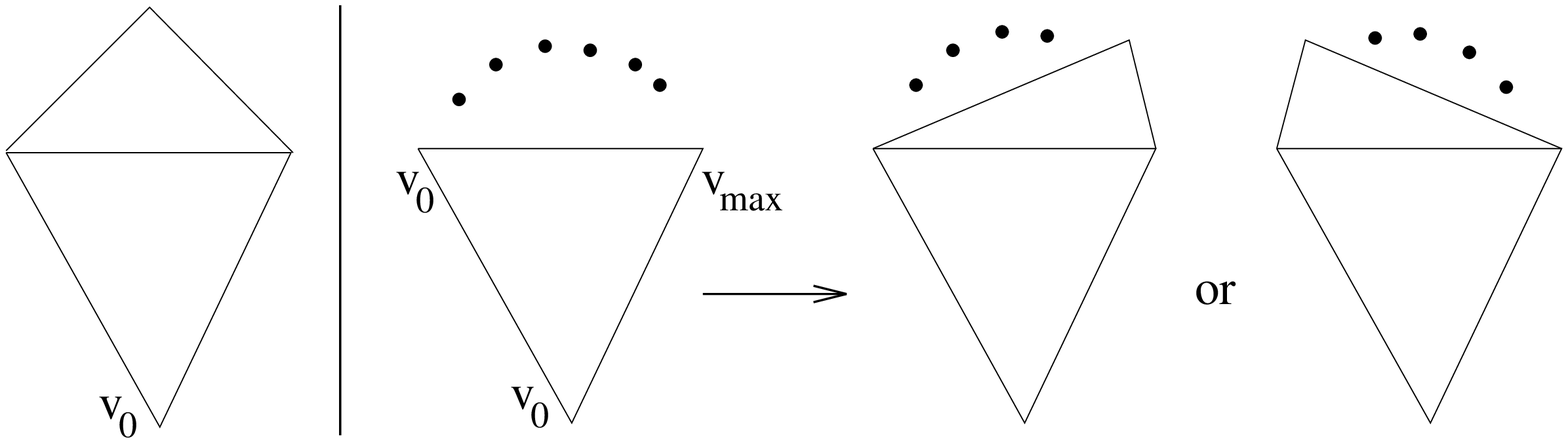}\hspace*{\fill}

\caption{Inductive definition of T-subdivisions}
\label{fig_tsubdivisions}
\end{figure} 

\section{P-resolutions and subdivisions}

We recall the inductive definition of graphs of type T using our
terminology (cf.~\cite{Stevens}). The graph $\Gamma (R(2,a,2))$ for
$a\geq 2$ is a graph of type T of the first kind. If $\Gamma (b)$ is a
graph of type T , then $\Gamma (\alpha b \gamma)$ and $\Gamma (\gamma
b \alpha)$ are graphs of type T of the second kind. We define a subset
$T$ of subdivions inductively as depicted in
Fig.\ref{fig_tsubdivisions}. If $(V,\Delta)\in T$, then there exists
exactly one vertex $v=v_i\in V$ such that
$\{v_{i-1},v_i,v_{i+1}\}\in\Delta$. This vertex is called the {\em
  central vertex} and $i$ the {\em central index} of $(V,\Delta)$.
By induction, one easily obtains:

\begin{lemma}
  Let $(V,\Delta)\in T$ and $(V,\Delta,a)\in \cS$ such that all
  vertices different from the central one have defect $0$. Then
  $\Gamma (R(a))$ is a graph of type T.
  
  On the other hand, assume that $\Gamma (R(a))$ is a graph of type T.
  Then there exists a subdivision $(V,\Delta)\in T$ such that
  $(V,\Delta,a)\in\cS$ and all vertices different from the central one
  have defect $0$. For each set $V$ of vertices with the appropriate
  cardinality there exists exactly one $\Delta$ with these properties.
\label{lem_tsubdivision}
\end{lemma}

\begin{lemma}
  Let $(V,\Delta,a)\in \cS$ with $(V,\Delta)\in T$ such that all
  vertices $v>v_0$ are saturated. Then
  there exist unique numbers $\ell,r\geq 0$ and a unique $a'\in\cW_2$
  such that $\Gamma (R(a'))$ is a graph of type T and $a=\alpha^\ell a'
  \alpha^r$.
  \label{lem_characterizeT}
\end{lemma}

A P-resolution of a quotient singularity $X$ is a partial resolution
$Y\longrightarrow X$ such that $K_Y\cdot E_i$ for all exceptional
divisors, and all singularities of $Y$ have a resolution graph of type
T or $A_k$ (cf. \cite{Stevens}, 3.1). An easy calculation with the
adjunction formula shows
\begin{lemma}
  Let $Y$ be a P-resolution of the cyclic quotient singularity $X$ and
  $G$ the graph of its minimal resolution $\Tilde{Y}$.  Let $J$ be the
  subset of nodes which blow down to singularities of type T in $Y$.
  Then each connected component of $J$ is a subgraph of $G$ being of
  type $T$. The $-1$-nodes in $G$ are adjacent to two nodes in $J$ and
  at least one of this node belongs to a component of $J$ which is of
  the second kind. The subset $K$ of nodes of $G$ which blow down to
  $A_k$--singularities can be characterized as follows: A $-2$-node
  belongs to $K$ if and only if none of its neighbours belongs to $J$.
  \label{lem_Presolution}
  
  On the other hand: If $G$ is a resolution graph of the quotient
  singularity $X$ and $J$ a subset of nodes of $G$ satisfying the
  above properties, then blowing down the nodes in $J$ and those in
  $K$ which is defined as above yields a P-resolution of $X$.
\end{lemma}

Assume $G=\Gamma (b)$ for a graph $G$ of the minimal resolution of a
P-resolution.  Let $I_1,\dots,I_r$ be the maximal intervals of $J$,
the set defined in lemma \ref{lem_Presolution}. Each interval $I_j$ is
associated with a interval $H_j$ of positions of $a=R(b)$.  These
intervals need not longer to be disjoint, they may have extremal
positions in common. Let $\{M_1,\dots,M_m\}$ be the set of all sets of
positions of $R(b)$ which are either equal to one of the $H_j$ or a
maximal interval in the complement of the union of all $H_j$. We
define inductively a subdivision $(V,\Delta)$ with
$(V,\Delta,a)\in\cS$. If $M_i$ is one of the $H_j$ then
$(V_i,\Delta_i)$ is chosen according to Lemma~\ref{lem_tsubdivision}
for the subgraph of $G$ determined by the interval $I_j$.  For all
other $M_i$ the subdivision $(V_i,\Delta_i)$ is chosen such that all
triangles have the minimal vertex as vertex. Then we set
\[(V,\Delta)=(V_1,\Delta)x \dots x(V_m,\Delta_m),\]
where $x$ denotes $\cdot$ or $*$, respectively, depending on whether
$M_i\cap M_{i+1} = \emptyset$ or not.

The set of all these labelled subdivisions is the set $\cP$.
Obviously, we can reconstruct the resolution graph $G$ and the subset
$J$ of nodes of $G$ from an element of $\cP$ which is associated
with $G$ by our construction, i.e.~there exists a canonical bijection
between P-resolutions and isomorphism classes of $\cP$.

On the other hand, we have the set $\cM$ of all labelled subdivisions
$(V,\Delta,a)$ with $a\in \cW_2$. The isomorphism classes of this set
are in a canonical one-to-one correspondence with graphs of minimal
resolutions of cyclic quotient singularities via the map $\Gamma\circ
R$ (cf.~\cite{Riemenschneider}).

The remaining part of this paper will provide a proof of:

\begin{theorem}
  There exists a canonical one-to-one correspondence between
  isomorphism classes of $\cM$ and isomorphism classes of $\cP$ with
  the following property. If $(V,\Delta,a)\in \cP$ and
  $(V',\Delta',a') \in \cM$ correspond to each other, then
  $\Gamma R(a')$ can be obtained from $\Gamma R(a)$ by
  successively blowing down $-1$-nodes.
\end{theorem}

\begin{proof}
  Blowing down and blowing up yields operations on the labels of
  subdivisions via the operator $R$, which is an involution. In order
  to find the right correspondence between $\cM$ and $\cP$, it is
  necessary to find suitable positions for blowing up and to define
  suitable modifications of the underlying subdivision. As we will see
  later, the subdivision itself gives a hint where to blow up.
  
  First of all, we prove a useful characterization of the set $\cP$:

\begin{lemma}
  A subdivision $(V,\Delta,a)$ is an element of $\cP$ if and only if
\begin{enumerate}
\item There exist no interior triangles.\label{lem_char1}
\item All vertices are saturated.\label{lem_char2}
\item If an irreducible component has more than two vertices then its
  label is in $\cW_2$.\label{lem_char3}
\item A vertex with label $1$ has two neighbouring vertices which
  belong to primary factors being elements of $T$. One of this factor
  has more than four vertices.\label{lem_char4}
\end{enumerate}\label{lem_characterize_cP}
\end{lemma}

\begin{proof}
  The necessity of these conditions follows easily from
  Lemma~\ref{lem_Presolution}. Let us assume that all these conditions
  hold.
  
  Condition~\ref{lem_char1}--\ref{lem_char3} imply with
  Lemma~\ref{lem_characterizeT} that a primary component
  $(V',\Delta')$ with more than two vertices is an element of $T$ and
  its label is $\alpha^\ell a' \alpha^r$ with $\Gamma R(a')$ a graph
  of type T and $\ell,r\geq 0$. If the next vertex to the lefthand or
  the righthand side of this component has label $1$ then $\ell = 0$
  or $r=0$, respectively, due to the definition of a saturated vertex.
  
  The conditons~\ref{lem_char2} and \ref{lem_char4} guarantee that the
  $-1$--nodes of $\Gamma (R(a))$ have the properties of lemma
  \ref{lem_Presolution}. Hence $\Gamma (R(a))$ is the graph of the
  minimal resolution of a P-resolution.
\end{proof}

This lemma shows that one should blow up in such a way that interior
triangles and non-saturated vertices disappear. In
Fig.~\ref{fig_blowup1} and Fig.~\ref{fig_blowup2}, one can see three
operations on subdivisions. We call the first one an operation of type
B1, and the second and third one operations of type B2a or B2b ,
respectively. All of them are operations of type B.

\begin{figure}

    \epsfysize4cm
    \hspace*{\fill}\epsffile{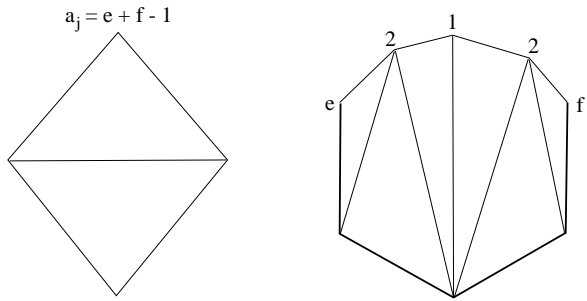}\hspace*{\fill}

\caption{An operation of type B1}
\label{fig_blowup1}
\end{figure}

\begin{figure}

    \epsfysize7.5cm
    \hspace*{\fill}\epsffile{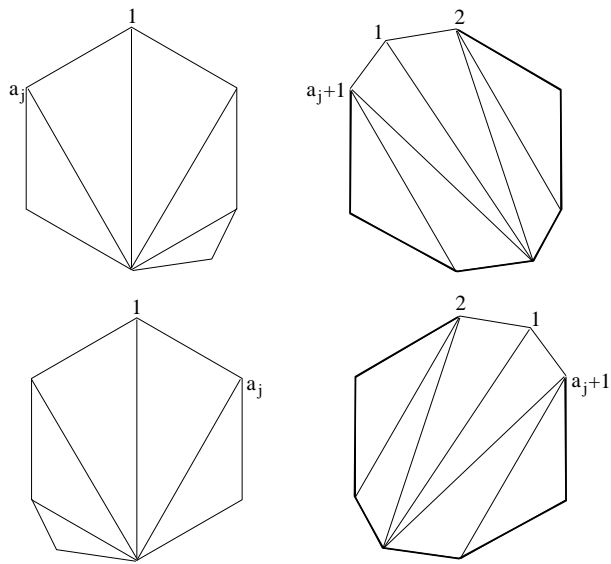}\hspace*{\fill}

\caption{Operations of type B2a and B2b}
\label{fig_blowup2}
\end{figure}

An operation of type B changes the label of the subdivision as
follows. From $a=a'(e+f-1)a''$ we proceed to $b=a'(e,2,1,2,f)a''$.
According to Lemma~\ref{lem_split}, we have $R(a)=R(a'(e))R((f)a'')$.
Therefore, $R(b) = R(a'(e+1))(1)R((f+1)a'')$, i.e.~$\Gamma (R(b))$ is
obtained by blowing up from $\Gamma (R(a))$. An analogous computation
shows that for operations of type B2a and B2b, the graph $\Gamma
(R(B))$ is obtained from the graph $\Gamma (R(a))$ by blowing up
between a $-1$-node and another node.

We will define a map $F:\cS \longrightarrow \cS$ which consists of
applying a suitable operation of type $B$, if appropriate, and show
that $F^k(V,\Delta,a)\in\cP$ for $(V,\Delta,a)\in\cM$ and $k$
sufficiently large. On the other hand, a map $G$ will be
defined with the property $F(G(V,\Delta,a))=(V,\Delta,a)
=G(F(V,\Delta,a))$ if $F(V,\Delta,a) \neq (V,\Delta,a)$ and
$G(V,\Delta,a)\neq (V,\Delta,a)$. Furthermore, $G^k
(V,\Delta,a)\in\cM$ for $(V,\Delta,a)\in\cP$ and $k$ sufficiently large.

Before we come to the definition of $F$ we need some technical
definitions. A configuration as depicted in Fig.~\ref{fig_hexafan} is
called a {\em hexagonal fan}. The vertex which is incident to all four
triangles is the {\em root} of the fan and the vertex opposite to the
root is the {\em apex} of the fan. Furthermore, we require that the
two triangles incident to the apex are not interior ones.

\begin{figure}

    \epsfysize2.5cm
\hspace*{\fill}\epsffile{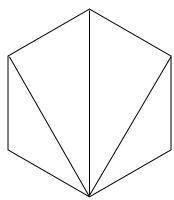}\hspace*{\fill}

\caption{A hexagonal fan}
\label{fig_hexafan}
\end{figure}

Assume  that we are given a subdivision $(V,\Delta)$ and a subset
$W$ of $V$ with $v_0\notin W$. The {\em height} of $W$ is the number
i+(k+1-j) where $v_i$ is the minimal element of $W$, $v_j$ the maximal
element of $W$ and $k+1=\# V$.

Let $(V,\Delta,a)\in\cS$. If
conditions~\ref{lem_char1}--\ref{lem_char3} hold, then
$F(V,\Delta,a):=(V,\Delta,a)$. Otherwise, there is a vertex $v_j$ such
that one of the following conditions hold.
\begin{enumerate}
\item There exists an interior triangle $\{v_i,v_j,v_k\}\in\Delta$
  with $i<j<k$.\label{case1}
\item The vertex $v_j$ is not saturated.\label{case2}
\item We have $a_j=1$ and $\{v_0,v_j\}$ is not an irreducible
  component of $(V,\Delta,a)$.\label{case3}
\end{enumerate}

Choose $j$ minimal with this property. Note that case~\ref{case1} and
case~\ref{case2} may occur simultaneously.  If condition~\ref{case1}
holds with $a_j\geq 3$ and all involved vertices have a label
different from $1$, then we apply an operation of type B1 to the
triangle $\{v_i,v_j,v_k\}$, otherwise we set $F(V,\Delta,a) =
(V,\Delta,a)$.  For this, we have to specify numbers $e,f$ with $a_j =
e + f -1$.  Choose $e$ such that it is the degree of the vertex whose
label it is in the new subdivision. After having applied this
operation either the number of interior triangles has become smaller
or the height of the first such triangle has become smaller.
Furthermore, we have obtained a hexagonal fan whose root is its
minimal or maximal vertex, its apex has label $1$ and all other
vertices of this fan have label different from $1$. Furthermore, no
interior triangle has this apex as vertex.

In the remaining cases condition~\ref{case2} holds in, we apply again
a operation of type B1. First we specify the triangle the operation is
applied to. Two cases can occur. First, that there exists a triangle
$\{v_{j-1},v_j,v_k\}$ with $k>j$, secondly, that there exists a
triangle $\{v_j, v_{j+1},v_k\}$ with $k<j$. In the first case, the
number $f$ is chosen such that it is the degree of the vertex whose
label it is in the new subdivision. In the second case, the number $e$
is chosen in such manner. We observe, that the number of non-saturated
vertices is reduced by this operation, since the vertex with
label $e$ or $f$, respectively, will be a vertex of a triangle with
three consecutive vertices, due to the construction.

We now deal with the case that condition~\ref{case3} holds.  We assume
that there is a hexagonal fan with $v_j$ as apex such that its root is
the minimal or maximal vertex of this fan and all vertices of the fan
being different from $v_j$ have a label $\geq 2$. In all other cases,
we define $F(V,\Delta,a)=(V,\Delta,a)$. If the root is minimal, we can
apply operation B2a, otherwise we apply operation B2b. Having applied
this operation, we have again a hexagonal fan as above. If the root is
not the minimal vertex $v_0$ then its height is smaller than the
height of the fan above.

If $(V,\Delta,a)$ is obtained from an element of $\cM$ by succesive
application of $F$, then proceeeding by induction, one can easily
verify the following facts:
\begin{itemize}
\item If $v_j$ is as above, then all vertices with label $1$ are
  smaller than this vertex.
\item Each vertex with label $1$ is apex of a hexagonal fan whose root
  is its minimal or maximal vertex. If $v_0$ is vertex of the
  hexagonal fan then it is its root.
\item A vertex with label $1$ has two neighbouring vertices which
  belong to primary factors being elements of $T$. One of these factors
  has more than four vertices.
\end{itemize}

Let us call a labelled subdivision with these properties a good one.
Hence, if $F(V,\Delta,a)=(V,\Delta,a)$ and $(V,\Delta,a)$ is good, then there exists no vertex
$v_j$ as above and hence $(V,\Delta,a)\in\cP$, as desired. The
considerations above show that this happens after finitely many
applications of $F$, since in each step we diminish one of the
following numbers: the number of interior triangles, the number of
non-saturated vertices, the number of vertices $v_j$ such that $a_j=1$
and $\{v_0,v_j\}$ is not an irreducible component, the height of the
first interior triangle, the height of the first hexafan whose apex
has label $1$.

Furthermore, the inverse map $G$ can be obtained as follows.  Apply on
good labelled subdivisions, if possible, the inverse of one of the
operations of type $B$ which yield a blowing down of the rightmost
$-1$-node in $\Gamma ( R (a))$. Otherwise, $G$ maps a labelled
subdivision onto itself. The desired properties of $G$ are now easily
to observe.
\end{proof}

\end{document}